\documentclass[]{article}
\usepackage{amssymb,amsmath}
\usepackage[mathscr]{euscript}
\input{epsf}

\newcounter{sec}

\newcounter{punct}[sec]

\def\punct{\refstepcounter{punct}{\arabic{sec}.\arabic{punct}.  }}

\newtheorem{theorem}{Theorem}[sec]

\newtheorem{lemma}[theorem]{Lemma}

\def\COUNTERS{\addtocounter{sec}{1}
              \setcounter{punct}{0}
          \setcounter{equation}{0}
          \setcounter{theorem}{0}
         }
          
          \def\sm{\smallskip}

\begin{document}

\newcommand{\supp}{\mathop {\mathrm {supp}}\nolimits}
\newcommand{\rk}{\mathop {\mathrm {rk}}\nolimits}
\newcommand{\Aut}{\mathop {\mathrm {Aut}}\nolimits}
\newcommand{\Out}{\mathop {\mathrm {Out}}\nolimits}
\newcommand{\OO}{\mathop {\mathrm {O}}\nolimits}
\renewcommand{\Re}{\mathop {\mathrm {Re}}\nolimits}
\newcommand{\ch}{\cosh}
\newcommand{\sh}{\sinh}

\def\0{\mathbf 0}

\def\ov{\overline}
\def\wh{\widehat}
\def\wt{\widetilde}

\renewcommand{\rk}{\mathop {\mathrm {rk}}\nolimits}
\renewcommand{\Aut}{\mathop {\mathrm {Aut}}\nolimits}
\renewcommand{\Re}{\mathop {\mathrm {Re}}\nolimits}
\renewcommand{\Im}{\mathop {\mathrm {Im}}\nolimits}
\newcommand{\sgn}{\mathop {\mathrm {sgn}}\nolimits}

\def\bfa{\mathbf a}
\def\bfb{\mathbf b}
\def\bfc{\mathbf c}
\def\bfd{\mathbf d}
\def\bfe{\mathbf e}
\def\bff{\mathbf f}
\def\bfg{\mathbf g}
\def\bfh{\mathbf h}
\def\bfi{\mathbf i}
\def\bfj{\mathbf j}
\def\bfk{\mathbf k}
\def\bfl{\mathbf l}
\def\bfm{\mathbf m}
\def\bfn{\mathbf n}
\def\bfo{\mathbf o}
\def\bfp{\mathbf p}
\def\bfq{\mathbf q}
\def\bfr{\mathbf r}
\def\bfs{\mathbf s}
\def\bft{\mathbf t}
\def\bfu{\mathbf u}
\def\bfv{\mathbf v}
\def\bfw{\mathbf w}
\def\bfx{\mathbf x}
\def\bfy{\mathbf y}
\def\bfz{\mathbf z}

\def\bfA{\mathbf A}
\def\bfB{\mathbf B}
\def\bfC{\mathbf C}
\def\bfD{\mathbf D}
\def\bfE{\mathbf E}
\def\bfF{\mathbf F}
\def\bfG{\mathbf G}
\def\bfH{\mathbf H}
\def\bfI{\mathbf I}
\def\bfJ{\mathbf J}
\def\bfK{\mathbf K}
\def\bfL{\mathbf L}
\def\bfM{\mathbf M}
\def\bfN{\mathbf N}
\def\bfO{\mathbf O}
\def\bfP{\mathbf P}
\def\bfQ{\mathbf Q}
\def\bfR{\mathbf R}
\def\bfS{\mathbf S}
\def\bfT{\mathbf T}
\def\bfU{\mathbf U}
\def\bfV{\mathbf V}
\def\bfW{\mathbf W}
\def\bfX{\mathbf X}
\def\bfY{\mathbf Y}
\def\bfZ{\mathbf Z}

\def\frD{\mathfrak D}
\def\frL{\mathfrak L}
\def\frG{\mathfrak G}
\def\frg{\mathfrak g}
\def\frh{\mathfrak h}
\def\frf{\mathfrak f}
\def\frl{\mathfrak l}

\def\bfw{\mathbf w}

\def\R {{\mathbb R }}
 \def\C {{\mathbb C }}
  \def\Z{{\mathbb Z}}
  \def\H{{\mathbb H}}
\def\K{{\mathbb K}}
\def\N{{\mathbb N}}
\def\Q{{\mathbb Q}}
\def\A{{\mathbb A}}

\def\T{\mathbb T}
\def\P{\mathbb P}

\def\G{\mathbb G}

\def\cD{\EuScript D}
\def\cL{\mathscr L}
\def\cK{\EuScript K}
\def\cM{\EuScript M}
\def\cN{\EuScript N}
\def\cP{\EuScript P}
\def\cQ{\EuScript Q}
\def\cR{\EuScript R}
\def\cW{\EuScript W}
\def\cY{\EuScript Y}
\def\cF{\EuScript F}
\def\cG{\EuScript G}
\def\cZ{\EuScript Z}
\def\cI{\EuScript I}
\def\cB{\EuScript B}
\def\cA{\EuScript A}
\def\cO{\EuScript O}

\def\bbA{\mathbb A}
\def\bbB{\mathbb B}
\def\bbD{\mathbb D}
\def\bbE{\mathbb E}
\def\bbF{\mathbb F}
\def\bbG{\mathbb G}
\def\bbI{\mathbb I}
\def\bbJ{\mathbb J}
\def\bbL{\mathbb L}
\def\bbM{\mathbb M}
\def\bbN{\mathbb N}
\def\bbO{\mathbb O}
\def\bbP{\mathbb P}
\def\bbQ{\mathbb Q}
\def\bbS{\mathbb S}
\def\bbT{\mathbb T}
\def\bbU{\mathbb U}
\def\bbV{\mathbb V}
\def\bbW{\mathbb W}
\def\bbX{\mathbb X}
\def\bbY{\mathbb Y}

\def\kappa{\varkappa}
\def\epsilon{\varepsilon}
\def\phi{\varphi}
\def\le{\leqslant}
\def\ge{\geqslant}

\def\B{\mathrm B}

\def\la{\langle}
\def\ra{\rangle}
\def\tri{\triangleright}

\def\lambdA{{\boldsymbol{\lambda}}}
\def\alphA{{\boldsymbol{\alpha}}}
\def\betA{{\boldsymbol{\beta}}}
\def\mU{{\boldsymbol{\mu}}}

\def\const{\mathrm{const}}
\def\rem{\mathrm{rem}}
\def\even{\mathrm{even}}
\def\SO{\mathrm{SO}}
\def\SL{\mathrm{SL}}
\def\SU{\mathrm{SU}}
\def\GL{\operatorname{GL}}
\def\End{\operatorname{End}}
\def\Mor{\operatorname{Mor}}
\def\Aut{\operatorname{Aut}}
\def\inv{\operatorname{inv}}
\def\red{\operatorname{red}}
\def\Ind{\operatorname{Ind}}
\def\dom{\operatorname{dom}}
\def\im{\operatorname{im}}
\def\md{\operatorname{mod\,}}

\def\ZZ{\mathbb{Z}_{p^\mu}}

\def\cH{\EuScript{H}}
\def\cQ{\EuScript{Q}}
\def\cL{\EuScript{L}}
\def\cX{\EuScript{X}}

\def\Di{\Diamond}
\def\di{\diamond}

\def\fin{\mathrm{fin}}
\def\ThetA{\boldsymbol {\Theta}}

\def\0{\boldsymbol{0}}

\def\F{\,{\vphantom{F}}_2F_1}
\def\FF{\,{\vphantom{F}}_3F_2}
\def\H{\,\vphantom{H}^{\phantom{\star}}_2 H_2^\star}
\def\HH{\,\vphantom{H}^{\phantom{\star}}_3 H_3^\star}
\def\Ho{\,\vphantom{H}_2 H_2}

\def\disc{\mathrm{disc}}
\def\cont{\mathrm{cont}}

\def\osigma{\ov\sigma}
\def\ot{\ov t}

\begin{center}

\Large

\bf
On Derkachov--Manashov $R$-matrices for\\ the principal series of unitary representations

\bigskip

\large
\sc
Yury A. Neretin%
\footnote{Supported by  FWF Der Wissenschaftsfonds (Austrian Scientific Fund), the grant P31591.}

\end{center}

\bigskip

{\small In 2001--2013 Derkachov and Manashov with coauthors  
 obtained  simple and natural expressions of  $R$-matrices
for the principal series of representations of the  groups $\mathrm{SL}(2,\mathbb{C})$,
 $\mathrm{SL}(2,\mathbb{R})$,  $\mathrm{SL}(n,\mathbb{C})$,   $\mathrm{SO}(1,n)$. The Yang--Baxter identities 
 for these intertwining operators are  kinds of  multivariate hypergeometric transformations.
 Derivations of the  identities are based on calculations
 'of physical level of rigor' with divergent integrals. 
 Our purpose  is a formal mathematical   justification of these results.}

\bigskip

\section{Introduction}

\COUNTERS

In several works starting 2001, Derkachov, Manashov, et al. (\cite{DKM}, \cite{KM}, \cite{DM}, \cite{CDS}, \cite{CDI}, \cite{VDI})  obtained natural constructions of $R$-matrices for the principal series of 
representations of the groups $\SL(2,\C)$,
 $\SL(2,\R)$,  $\SL(n,\C)$ and the spherical principal series of $\SO(1,n)$. Derivations were given on 'physical level of rigor', and these formulas remain to be a topic of 'mathematical physics'.
 In author's opinion, these objects are  important for non-commutative harmonic analysis, it is sad that they were 'lost' for this branch of mathematics. 
The purpose of this note is  to justify formally the derivation of the Yang-Baxter identity. We partially achieve this purpose and indicate  obstacles for an obtaining of a perfect result. In author's opinion, a further justification is an
interesting problem for functional analysis.


\sm

{\bf\punct The principal series of unitary representations of $\SL(2,\R)$.}
See \cite{GGV}, Chapter 7.
Denote by $\T$ the circle $|z|=1$ on the complex plane.
Denote by $\SU(1,1)$ the group of all complex matrices $g=\begin{pmatrix}
                                                           a&b\\\ov b &\ov a
                                                          \end{pmatrix}$
such that $|a|^2-|b|^2=1$. Recall that this group is isomorphic to $\SL(2,\R)$.

The group $\SU(1,1)$  acts on $\T$ by the M\"obius transformations 
$$
z\mapsto \frac{b+z\ov a}{a+z\ov b}
.$$
For $p\in \C$ we define a {\it representation $T_p(g)$ of the principal series}%
\footnote{More precisely, this is the  {\it even principal series} since our representations  are trivial on the center of $\SU(1,1)$ consisting 
of matrices $\begin{pmatrix}\pm1 &0\\0& \pm1 \end{pmatrix}$.}
of $\SU(1,1)$ in the space $C^\infty(\T)$ of $C^\infty$-smooth functions
on $\T$ by the formula
\begin{equation}
T_p\begin{pmatrix} a&b\\\ov b &\ov a\end{pmatrix}f(z)=
f\Bigl(\frac{b+z\ov a}{a+z\ov b} \Bigr)    \bigl|a+z\ov b\bigr|^{-1+p}.
\label{eq:Tp-sl}
\end{equation}
For $p\in i\R$ these representations are unitary in $L^2(\T)$ and they are called {\it representations
of the unitary principal series}.

For $p\in \C$ we define the operator $J(p):C^\infty(\T)\to C^\infty(\T)$ by
$$
J(p) f(z)= \frac{1}{2\,C(p)}\int_{|u|=1} |z-u|^{-1+p}  f(u)\,\frac{du}{iu},
$$
where 
\begin{equation}
C(p):=\Gamma(p)\cos\pi p/2.
\label{eq:C}
\end{equation}
Notice that  $\frac{du}{iu}$ is the Lebesgue measure on $\T$. 

\sm

 A straightforward calculation shows that
{\it the operator $J(p)$ intertwines representations $T_p$ and $T_{-p}$},
$$
T_p(g)\, J(p)=J(p) \, T_{-p}(g).
$$

\sm

The operators $J(p)$ satisfy the following properties:

\sm

1) The functions $z^n$ are eigenfunctions of $J(p)$,
\begin{equation}
 J(p):\,\, z^n\mapsto \lambda_n(p)\,z^n,\quad\text{where} \quad \lambda_n(p):=2^{-i\Im p}\, \frac{\Gamma(1/2-p-n)}{\Gamma(1/2-p+n)}
 \label{eq:eigenfunctions}
\end{equation}
(this can be derived by elementary calculations from \cite{PBM1}, formulas (2.5.12.11), (2.5.12.36), or \cite{HTF}, formula (1.5.1.29)).

\sm

2) This implies that $p\mapsto J(p)$ is a meromorphic operator-valued function. If $p\ne\pm 1$, $\pm3$, $\pm5$, \dots,
then $J(p)$ is an isomorphism of representations in $C^\infty(\T)$.

\sm

3) Also, we see that  for $p\in i\R$ this operator
is unitary in $L^2(\T)$.   


\sm

4) Applying  asymptotic formula \cite{HTF}, (1.18.4), we get
$$
\lambda_k(p)= 2^{-i\Im p}\,
|k|^{-\Re p} 
\bigl(1+O(|k|^{-1})\bigr), \qquad \text{as $k\to\pm\infty$},
$$
where $O(\cdot)$ is locally uniform in $p$.

5) So, for $\Re p\ge 0$, $p\ne 1$, 3, 5, $\dots$ the operator $J(p)$ is bounded and compact in $L^2(\T)$.
This operator-valued function is continuous in the weak operator topology up to the line $\Re p=0$.


\sm

{\bf \punct Symmetries of tensor products.} See, \cite{Ner}, \cite{NO}, Sect. 7.
A function
$|z-u|^\alpha$ on the torus $\T^2=\T\times \T$ satisfies the identity
\begin{equation}
\Bigl| \frac{\ov b+z a}{a+z\ov b}-\frac{\ov b+u a}{a+u\ov b}\Bigr|^\alpha=|z-u|^\alpha \,|a+z\ov b|^{-\alpha} \,|a+u\ov b|^{-\alpha}.
\label{eq:z-u}
\end{equation}

Consider the following operator $A(\alpha)$ in a space of functions on torus $\T^2$:
$$
A(\alpha) f(z,u)=|z-u|^\alpha\,f(z,u).
$$
For $\alpha\in \R$ the  operator $A(\alpha)$ is unitary in $L^2(\T^2)$. 
For $\alpha\notin \N$ this operator is not well-defined in $C^\infty(\T^2)$.
Let us forget about analytic difficulties
related to a definition of $A(\alpha)$.
Formally, by \eqref{eq:z-u},
$$
 \bigl(T_p(g)\otimes T_q(g)\bigr)\, A(\alpha)=A(\alpha)\, \bigl(T_{p-\alpha}(g)\otimes T_{q-\alpha}\bigr). 
$$
We also have
$$
 \bigl(T_p(g)\otimes T_q(g)\bigr)\cdot \bigl(J(p)\otimes J(q)\bigr)= \bigl(J(p)\otimes J(q)\bigr) \cdot  \bigl(T_{-p}(g)\otimes T_{-q}(g)\bigr).
$$

\begin{figure}
 $$\epsfbox{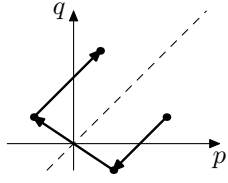}$$
 \caption{To the definition of $R$-matrix.}
 \label{fig:de}
\end{figure}

{\bf\punct The  $R$-matrix.}
For a complex $\sigma$
consider the operator  $\check R(\sigma)$ in a space of functions on the torus $\T^2$
given by
\begin{multline}
\check R(\sigma)\,f(y_1,y_2)=\check R^{p,q}(\sigma)f(y_1,y_2):=\\
\frac1{4\, C\bigl(-\sigma+\tfrac12(q-p)\bigr)\, C\bigl(-\sigma+\tfrac12(p-q)\bigr)}
\times\\ \times
\int\limits_{\T^2}
|y_1-y_2|^{\sigma+\frac12(p+q)}|x_1-y_1|^{-1-\sigma+\frac12(q-p)}|x_2-y_2|^{-1-\sigma+\frac12(p-q)}
|x_1-x_2|^{\sigma-\frac12(p+q)} \times\\ \times f(x_1,x_2)\,\frac{dx_1}{ix_1}\,\frac{dx_2}{ix_2}=
\label{eq:R}
\end{multline}
\begin{multline*}=
[\text{Gamma-factor}]\cdot |y_1-y_2|^{\sigma+\frac12(p+q)}
 \times\\ \times
 \int\limits_{\T^2}
|x_1-y_1|^{-1-\sigma+\frac12(q-p)}\,|x_2-y_2|^{-1-\sigma+\frac12(p-q)}
 \times\\ \times
\Bigl(|x_1-x_2|^{\sigma-\frac12(p+q)} \, f(x_1,x_2)\Bigr)\,\frac{dx_1}{ix_1}\,\frac{dx_2}{ix_2},
\end{multline*}
or, equivalently,
\begin{multline}
\check R(\sigma)= A\bigl( \sigma+\tfrac12(p+q)\bigr)\, \Bigl( J\bigl(-\sigma+\tfrac12(q-p)\bigr)\otimes
 J\bigl(-\sigma+\tfrac12(p-q)\bigr)\Bigr)
 \times\\ \times
 A\bigl( \sigma-\tfrac12(p+q)\bigr).
\label{eq:RR}
\end{multline}
Then formally  (see Fig.\ref{fig:de})%
\footnote{It seems that $\check R(\sigma)$ is  the simplest intertwining operator 
$T_p\otimes T_q\to T_q\otimes T_p$ different from the permutation of tensor factors.}
$$
\check R (\sigma) \, \bigl(T_p(g)\otimes T_q(g)\bigr)= \bigl(T_q(g)\otimes T_p(g)\bigr)\, \check R(\sigma).
$$

Next, consider a representation $T_p\otimes T_q\otimes T_r$ in a space of functions 
on 3-dimensional torus $\T^3$.
For $k$, $l=1$, 2, 3, where $k\ne l$, denote by $\check R_{kl}(\sigma)$ the operator $\check R(\sigma)$ acting in $k$-th and $l$-th factors.

Kirch and Manashov \cite{DM} established on  physical level of rigor the following {\it Yang--Baxter identity}:
\begin{equation}
\check  R_{23}(\theta-\tau)\, \check R_{12}(\theta)\, \check R_{23}(\tau)=
\check R_{12}(\tau)\,\check R_{23}(\theta)\,\check R_{12}(\theta-\tau),
 \label{eq:RRR}
\end{equation}
or (with upper superscripts)
$$
\check  R_{23}^{p,q}(\theta-\tau)\,\check  R_{12}^{p,r}(\theta)\,\check  R_{23}^{q,r}(\tau)=
\check  R_{12}^{q,r}(\tau)\,\check R_{23}^{p,r}(\theta)\,\check R_{12}^{p,q}(\theta-\tau).
$$
A similar construction remains valid for the group $\SL(2,\C)$  and, more generally, for the groups $\SL(n,\C)$
 (Derkachov, Manashov \cite{DM}).
 Finite-dimensional representations of these groups are realized at certain points of degenerations of the  principal series;
according \cite{CDS} for finite-dimensional representations this construction produces classical $R$-matrices.

A translation of statements about operators $\check R(\sigma)$  to 'mathematical level of rigor' is not self-evident, 
since the operators $\check R(\sigma)$ are very singular. They are not well-defined in standard functional spaces 
of Sobolev type,
and it is not clear how to define their products.

In Section 2 we prove a weaker statement:

\begin{theorem}
	\label{th:1}
 Let $T_p$, $T_q$, $T_r$ be representations of the unitary principal series,
 i.~e., $p$, $q$, $r\in i\R$. Let $\theta$, $\tau\in i\R$. Then the identity \eqref{eq:RRR} holds.
\end{theorem}

For a proof we justify the star-triangle identity from \cite{KM}, \cite{DM} for unitary operators (see Theorem \ref{th:star-triangle}).

In Section 3, for complex values of parameters we show that 
 all operators and all subexpressions, which are present in formula \eqref{eq:RRR}, are well-defined
 as operators sending smooth functions to distributions. However, we can not define products of such operators.

\section{Proof of the Yang-Baxter property}

\COUNTERS

{\bf \punct A beta-integral.} See \cite{KM}.

\begin{lemma}
 Let $\Re \alpha$, $\Re\beta$, $\Re \gamma >0$, $\alpha+\beta+\gamma=1$. Let $a$, $b$, $c\in \T$.
 Then
 \begin{multline}
 \int_{|z|=1} |z-a|^{\alpha-1}\, |z-b|^{\beta-1}\, |z-c|^{\gamma-1}\,\frac{dz}{iz}=
 \\= \frac4\pi C(\alpha)\, C(\beta)\, C(\gamma)\,
 |a-b|^{-\gamma}\, |b-c|^{-\alpha} \,|a-c|^{-\beta}
 =\\=
 \frac{2\,C(\alpha)\,C(\beta)}{C(\alpha+\beta)}\, |a-b|^{-\gamma}\, |b-c|^{-\alpha}\, |a-c|^{-\beta}
 , 
 \label{eq:beta}
 \end{multline}
 where
 $
 C(\delta):=\Gamma(\delta) \cos \pi\delta/2
 $ is the same as above {\rm(\ref{eq:C})}.
\end{lemma}

{\sc Proof.}
First,
\begin{multline*}
 \int_\R |t|^{\mu-1} |1-t|^{\nu-1}dt=
 \int_{-\infty}^0+\int_0^1+\int_1^\infty
 =\\=
 \frac{\Gamma(\mu)\,\Gamma(\nu)}{\Gamma(\mu+\nu)}+
 \frac{\Gamma(\nu)\,\Gamma(1-\mu-\nu)}{\Gamma(1-\mu)}+
 \frac{\Gamma(\mu)\,\Gamma(1-\mu-\nu)}{\Gamma(1-\nu)}=\\=
 \frac{\Gamma(\mu)\Gamma(\nu)}{\Gamma(\mu+\nu)} \frac{2\cos\pi\mu/2\, \cos\pi\nu/2}{\cos\pi(\mu+\nu)/2}
 =\\=
\frac r\pi \Gamma(\mu)\,\Gamma(\nu)\Gamma(1-\mu-\nu)\, \cos\pi\mu/2\, \cos\pi\nu/2\, \cos\pi(1-\mu-\nu)/2
\end{multline*}
(we used the standard beta-integral and the beta-integral \cite{PBM1}, 2.2.4.24).

Second, in (\ref{eq:beta})
we set
$$
z=\frac{a(c-b)+t c(b-a)}{(c-b)+t(b-a)}
$$
and reduce this integral to (\ref{eq:beta}).
\hfill $\square$

\sm

{\bf \punct The star-triangle identity.}
Consider a triple tensor product
$T_p\otimes T_q\otimes T_r$ realized in the space of functions
of variables $(x_1,x_2,x_3)\in \T^3$.

For $k\ne l$ consider the operators
$$
A_{k,l}(\alpha)\,f(x_1,x_2,x_3)
:=|x_k-x_l|^\alpha\, f(x_1,x_2,x_3).
$$
For each $k$ we define the operator $J_k(\beta)$ as the operator
$J(\beta)$ acting on the coordinate $x_k$. For $\alpha$, $\beta\in i\R$
these operators are unitary (see \eqref{eq:eigenfunctions}).

In this section, we consider $p$, $q$, $r\in i\R$ and $\alpha$, $\beta$, \dots $\in i\R$. 
We have the following obvious commutativity relations%
\footnote{The identity \eqref{eq:JJ} actually is used below only for $k\ne l$.}:
\begin{align}
 A_{k,l}(\alpha)\, A_{k',l'}(\beta)&=A_{k',l'}(\beta) \, A_{k,l}(\alpha);
 \label{eq:comm1}
 \\
 J_k(\alpha)\,J_l(\beta)&=J_l(\beta)\, J_k(\alpha);
 \label{eq:JJ}
 \\
 A_{k,l}(\alpha)\, J_m(\beta)&=  J_m(\beta)\, A_{k,l}(\alpha),\quad \text{where $k$, $l$, $m$ are distinct.}
  \label{eq:comm3}
\end{align}

\begin{theorem}
\label{th:star-triangle}
 Let parameters $\alpha$, $\beta$, $\gamma$ range in the  set 
\begin{equation} 
\Sigma:\, \text{$\alpha$, $\beta$, $\gamma\in i\R$ and $\alpha+\beta+\gamma=0$.}
\label{eq:domain-Delta}
 \end{equation}
 Then the following {\rm star-triangle identity} for intertwining operators holds 
 \begin{equation}
  J_1(\alpha)\, A_{12}(\beta)\, J_1(\gamma)=
  A_{12}(-\gamma)\, J_1(-\beta)\, A_{12} (-\alpha).
  \label{eq:star-triangle}
 \end{equation}
\end{theorem}

{\sc Proof.} 
 Denote by $C^\infty_c(\T^2\setminus\Delta)$ the space of smooth functions on
the torus $\T^2$, which vanish in a neighbourhood of  the diagonal $x_1=x_2$ of $\T^2$.
A sequence $\Phi_j\in C^\infty_c(\T^2\setminus \Delta)$ converges to $\Phi$
if $\Phi_j$ uniformly converges to $\Phi$ with all partial derivatives and there is a neighbourhood
$\cO$ of $\Delta$, where all functions $\Phi_j$ are zero.
The  space $C^\infty_c(\T^2\setminus \Delta)$ is dense in $L^2(\T)$.

Recall that for any operator $Q$ in a Hilbert space we have an Hermitian form
$\la Q\Psi,   \Phi \ra$. Let us verify that such forms coincide for both sides of 
\eqref{eq:star-triangle},
 \begin{multline}
 \bigl\la J_1(\alpha)\, A_{12}(\beta)\, J_1(\gamma)\Psi,\,\Phi\bigr\ra_{L^2(\T^2)}
=\\=
\bigl\la   A_{12}(-\gamma)\, J_1(-\beta)\, A_{12} (-\alpha)\Psi,\,\Phi\bigr\ra_{L^2(\T^2)}.
  \label{eq:star-triangle-0}
 \end{multline}
Since operators are unitary, it is sufficient to check this for functions 
$\Phi$, $\Psi\in C^\infty_c(\T^2\setminus\Delta)$.
We will verify the identity representing it in the form
\begin{equation}
 \la  A_{12}(\beta) \, J_1(\gamma)\Psi,\, J_1(\alpha)^* \Phi \ra=
\la   J_1(-\beta)\, A_{12} (-\alpha)\Psi,\, A_{12}(-\gamma)^*\Phi \ra.
 \label{eq:star-triangle-1}
\end{equation}

Consider a wider domain for our parameters:
\begin{equation}
\Xi:\,\Re \alpha\ge 0, \quad -1/2<\Re\beta\le 0,\quad \Re \gamma\ge0, \quad
\alpha+\beta+\gamma=0
\label{eq:domain-Xi}
\end{equation}
(in particular, $\Re\alpha$, $\Re \gamma < 1/2$).
Consider its interior
\begin{equation}
\Xi^\circ:\,\Re \alpha> 0, \quad -1/2<\Re\beta< 0,\quad \Re \gamma>0, \quad
\alpha+\beta+\gamma=0.
\label{eq:domain-Xi-circ}
\end{equation}
The domain $\Sigma$, see \eqref{eq:domain-Delta}, interesting for us
is contained in the boundary
$\Xi\setminus \Xi^\circ$.

Our theorem is a corollary of the following two lemmas.

\begin{lemma}
\label{l:continuity-1}
For fixed $\Phi$, $\Psi\in C^\infty_c(\T^2\setminus\Delta)$ and for $\alpha$, $\beta$, $\gamma$ ranging in the domain $\Xi$
the inner products in both sides of {\rm (\ref{eq:star-triangle-1})} are continuous functions 
of parameters $\alpha$, $\beta$, $\gamma$. 
\end{lemma}

\begin{lemma}
\label{l:continuity-2}
For fixed $\Phi$, $\Psi\in C^\infty_c(\T^2)$ and for $\alpha$, $\beta$, $\gamma$ ranging in the domain $\Xi^\circ$
the inner products in both sides of {\rm(\ref{eq:star-triangle-1})} coincide.
\end{lemma}

{\sc Proof of Lemma \ref{l:continuity-1}.}
Below we use the following statement, see 
\cite{Rud}, Theorem 2.17:

\sm

--- {\it  Let $X$ be a Fr\'echet space {\rm(}in particular, a Banach space{\rm)}, $Y$, $Z$ be  topological linear spaces.
Let $F:X\times Y\to Z$ be a separately continuous bilinear map. Then $F$ is jointly sequentially continuous, 
i.e.,  
$$\text{$x_n\to x$ in $X$ and $y_n\to y$ in $Y$ imply $F(x_n, y_n)\to F(x,y)$ in $Z$.}
$$
If $Y$ is metrizable, then $F$ is jointly continuous.
}

\sm

 We formulate the following chain of remarks.

\sm 

--- The maps 
$$
\Psi\mapsto A_{12} (-\alpha)\Psi,\quad\Phi\mapsto A_{12}(-\gamma)^*\Phi
$$
are well-defined maps $C^\infty_c(\T^2\setminus\Delta)\to C^\infty_c(\T^2\setminus \Delta)$
and images of $\Phi$, $\Psi$ continuously depend  on parameters $\alpha$, $\gamma$.

\sm 

--- The maps 
$$
\Phi\mapsto J_1(\alpha)^* \Phi, \quad \Psi\mapsto J_1(\gamma)\Psi
$$
are well-defined maps $C^\infty(\T^2)\to C^\infty(\T^2)$ and functions
  $J_1(\alpha)^* \Phi$, $J_1(\gamma)\Psi\in  C^\infty(\T^2)$ continuously  depend on parameters $\alpha$, $\gamma$.
Indeed, the  operator $J_1(\alpha)$ is a  convolution with a well-defined 
distribution 
$$(\Gamma(\alpha+1)\cos(\alpha\pi))^{-1}(1-x_1)^{-1+\alpha}$$
 on the circle meromorphically depending on $\alpha$
(see, e.g., \cite{GSh}, Sect. 1.3, 1.5). Also, we can refer to formula \eqref{eq:eigenfunctions}.

\sm 

--- For the same reason, the map
$$
\Psi\mapsto J_1(-\beta) A_{12} (-\alpha)\Psi
$$
is a well-defined map $C^\infty_c(\T^2\setminus\Delta)\to C^\infty(\T^2)$ and the image of $\Psi$ continuously depends on parameters $\beta$, $\alpha$.
More precisely, $A_{12} (-\alpha)\Psi$ is a continuous family in a Hilbert space $L^2(\T^2)$,
and $J_1(-\beta)$ is a continuous family in the (non-metrizable) space of bounded operators
$L^2(\T^2)\to L^2(\T^2)$ equipped with the weak operator topology. So we can use the statement about joint sequential continuity.

\sm 

--- The map
$$\Psi\mapsto A_{12}(\beta) J_1(\gamma)\Psi$$
is a well-defined map $C^\infty(\T^2)\to L^2(\T^2)$ (since $\Re\beta>-1/2$) and the image of $\Psi$
continuously depend on parameters $\beta$, $\gamma$. More precisely, 
$J_1(\gamma)\Psi$ is a continuous family in the Banach space $C(\T^2)$ of continuous functions,   and
$A_{12}(\beta)$ is a continuous family   
  in the Banach space of bounded operators $C(\T^2)\to L^2(\T^2)$. So we have the joint continuity. 

\sm 

--- Finally, for fixed functions $\Phi$, $\Psi\in C^\infty_c(\T^2\setminus \Delta)$
the expressions
\begin{equation*}
 J_1(\alpha)^* \Phi,\qquad A_{12}(\beta) J_1(\gamma)\Psi, \qquad
  A_{12}(-\gamma)^*\Phi, \qquad J_1(-\beta) A_{12} (-\alpha)\Psi
\end{equation*}
considered as elements of $L^2(\T^2)$ continuously depend on parameters $\alpha$, $\beta$, $\gamma$.

\sm

Thus, the inner products in both sides of (\ref{eq:star-triangle-1}) are continuous functions 
of parameters $\alpha$, $\beta$, $\gamma$. \hfill $\square$

\sm 

Consider the left hand side of the equality \eqref{eq:star-triangle-0} (or, equivalently,
 \eqref{eq:star-triangle-1}), i.e., the iterated integral:
\begin{multline}
\frac{1}{4C(\alpha)C(\gamma)}\times\\
 \int_{\T^4} \ov {\Phi(\eta,x_2)} |\xi-\eta|^{-1+\alpha}|\xi-x_2|^{\beta} |\xi-x_1|^{-1+\gamma}
 \Psi(x_1,x_2)\frac{dx_1}{ix_1}\frac{d\xi}{i\xi}\frac{d\eta}{i\eta}\frac{d x_2}{i x_2}.
 \label{eq:int-Phi-Psi}
\end{multline}

\begin{lemma}
\label{l:convergence}
For $(\alpha,\beta,\gamma)\in \Xi^\circ$, the integral \eqref{eq:int-Phi-Psi}
absolutely converges.
\end{lemma}

{\sc Proof of Lemma \ref{l:convergence}.} We must verify the convergence of the integral
$$
 \int_{\T^4}  |\xi-\eta|^{-1+\Re \alpha}\,|\xi-x_2|^{\Re \beta} \,|\xi-x_1|^{-1+\Re\gamma}\,
 \frac{dx_1}{ix_1}\frac{d\xi}{i\xi}\frac{d\eta}{i\eta}\frac{d x_2}{i x_2}.
$$
Integrating in $\xi$ using the beta-integral \eqref{eq:beta}, we get
$$
\mathrm{const}(\alpha,\beta,\gamma)\cdot \int_{\T^3} 
|x_2-\eta|^{-\Re\gamma} \,|x_1-\eta|^{-1-\Re\beta}\,|x_1-x_2|^{-\Re\alpha}\,
 \frac{dx_1}{ix_1}\frac{d\eta}{i\eta}\frac{d x_2}{i x_2}.
$$
We come to a question about integrability  of the following function of two real 
variables
$$
|u|^a\, |v|^b\, |u-v|^c
$$
near zero. Passing to polar coordinates we  get sufficient and necessary conditions: $\Re a$, $\Re b$, $\Re c>-1$ and $\Re(a+b+c)>-2$,
or, in the  notation of this proof,
$$
\Re\alpha<1;\quad
\Re\beta<0;\quad
\Re \gamma<1,\quad
\Re(\alpha+\beta+\gamma)<1.
$$
This is valid in the domain $\Xi$.
\hfill $\square$

\sm

{\sc Proof of Lemma \ref{l:continuity-2}.}
We  integrate in $\xi$ and get
\begin{multline*}
\frac{1}{4C(\alpha)C(\gamma)}  \frac{2(\alpha)C(\gamma)}{C(\beta)}\times\\\times
 \int_{\T^4} \ov {\Phi(\eta,x_2)} \,
|x_2-\eta|^{-\gamma} |x_1-\eta|^{-1-\beta}|x_1-x_2|^{-\alpha}
 \Psi(x_1,x_2)\,\frac{dx_1}{ix_1}\frac{d\eta}{i\eta}\frac{d x_2}{i x_2},
\end{multline*}
i.e., the right hand side of \eqref{eq:star-triangle-0}.
\hfill $\square$




\sm

{\bf\punct Proof of Theorem \ref{th:1}.}
In the left hand side of (\ref{eq:RRR}) we have the product of the following  operators:
\begin{align*}
 \check R_{23}(\tau):T_p\otimes T_q \otimes T_r\to T_p\otimes T_r \otimes T_q;
 \\
 \check R_{12}(\theta): T_p\otimes T_r \otimes T_q\to T_r\otimes T_p \otimes T_q;
 \\
 \check R_{23}(\theta-\tau):T_r\otimes T_p \otimes T_q\to T_r\otimes T_q \otimes T_p.
\end{align*}

Keeping in the mind \eqref{eq:eigenfunctions}, we present the left hand side of the Yang--Baxter equation (\ref{eq:RRR})
 in the following form:
\begin{multline*}
 \check R_{23}(\theta-\tau)\, \check R_{12}(\theta)\, \check R_{23}(\tau)=
 \\
 A_{23}\bigl(\theta-\tau+\tfrac12(p+q)\bigr)J_2\bigl(-\theta+\tau+\tfrac12(q-p)\bigr)J_3\bigl(-\theta+\tau+\tfrac12(p-q)\bigr)
 A_{23}\bigl(\theta-\tau-\tfrac12(p+q)\bigr)\times\\
 \times
 A_{12}\bigl(\theta+\tfrac12(p+r)\bigr) J_1\bigl(-\theta+\tfrac12(r-p)\bigr)J_2(-\theta +\tfrac12(p-r)\bigr) A_{12}(\theta-\tfrac12(p+r)\bigr)
 \times\\
 \times A_{23}\bigl(\tau+\tfrac12(q+r)\bigr)J_2\bigl(-\tau+\tfrac12(r-q) \bigr)J_3\bigl(-\tau+\tfrac12(q-r)\bigr) A_{23}\bigl(\tau-\tfrac12(q+r)\bigr).
 \label{eq:YB}
\end{multline*}
Using commutativity relations \eqref{eq:comm1}--\eqref{eq:comm3} and applying 8 times the star-triangle identity,
we can get the right hand side of the Yang--Baxter equation. It is difficult to understand a chain of
15 such formulas. In fact,  we  confirm
the  diagram calculation from \cite{KM}, \cite{DM},
  transformations of divergent integrals 
can be understood as transformations of products of unitary operators (but we must additionally indicate order of factors on diagrams, see below; in my opinion this is necessary).

\begin{figure}
 $$J_2(\alpha): \quad\epsfbox{fey4.1}
 \qquad\quad A_{23}(\beta):\quad \epsfbox{fey4.2}
 $$
 \caption{Notation for Fig. \ref{fig:star-triangle}-\ref{fig:predposl}. Operators and corresponding elements of diagrams. 
For $J(\cdot)$ labels on the ends of the edge coincide, for $A(\cdot)$ they are different.} 
 \label{fig:add1}
\end{figure}

\begin{figure}
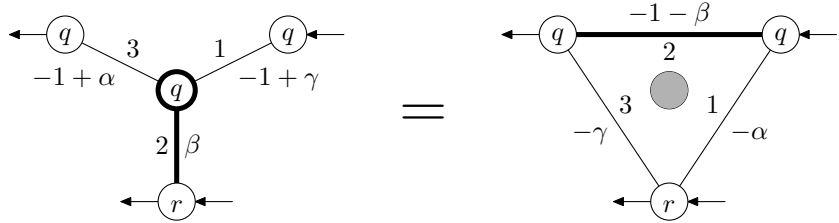

 $$ \epsfbox{fey5.2} \qquad
\begin{matrix}
 \text{\Huge $=$}\\
 \phantom{\biggl|}\\
  \phantom{\biggl|}\\
   \phantom{\biggl|}
\end{matrix}
\qquad \epsfbox{fey5.1} $$
\caption{Notation for Fig. \ref{fig:1}-\ref{fig:predposl}. The star-triangle relation, $\boxed{\alpha+\beta+\gamma=0}$.
In the left hand side ('{\it star}') we emphasize the central vertex in bold (to visualize a place of transformation) and also emphasize 
the edge, where labels at vertices are different. In the right hand side we emphasize the edge, where
labels of vertices coincide and mark transformed '{\it triangle}' by a grey circle.
Our operators act from one copy $L^2(\T^2)$ to another copy of $L^2(\T^2)$. We mark variables of
the initial torus by arrows entering into circles, variables of the target torus by outgoing arrows.
}
\label{fig:star-triangle}
\end{figure}

\begin{figure}
$$\epsfbox{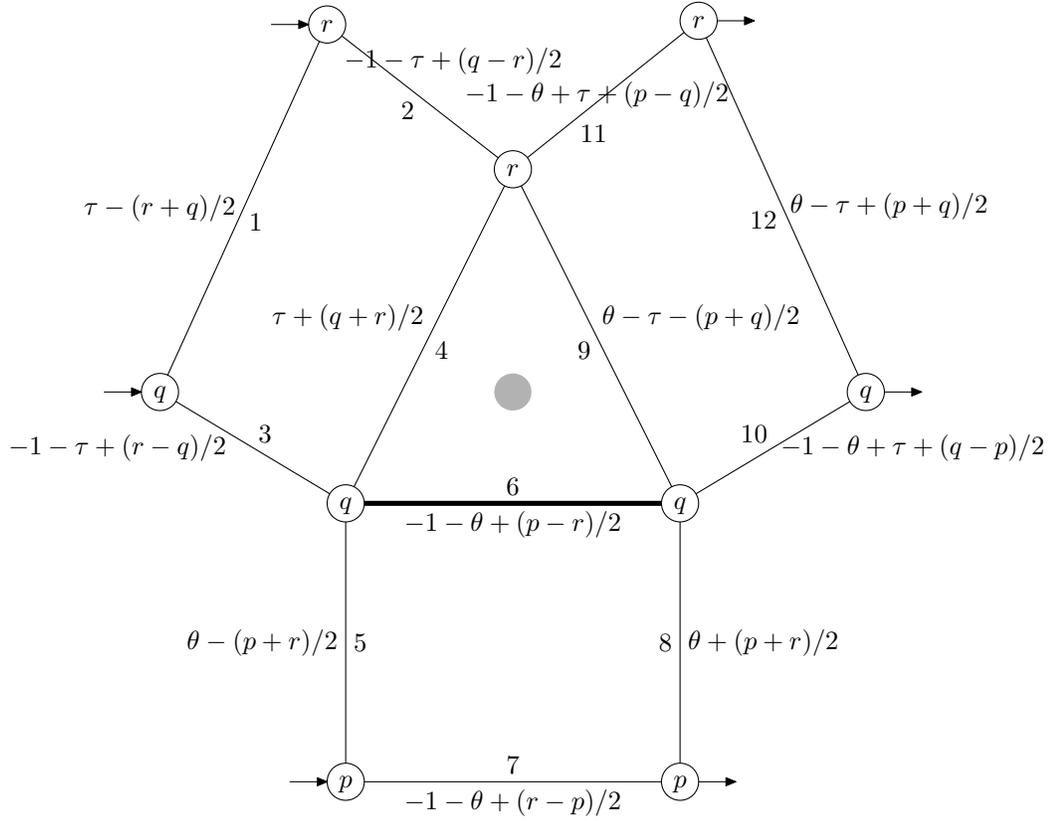}$$
\caption{The diagram corresponding to $\check R_{23}(\theta-\tau) \check R_{12}(\theta) \check R_{23}(\tau)$.}
\label{fig:1}
\end{figure}

\begin{figure}
$$\epsfbox{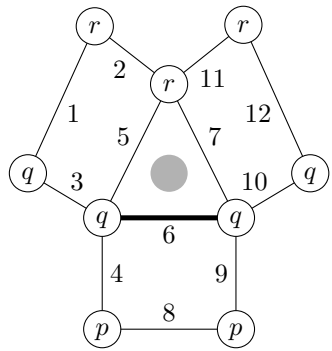}$$
\caption{After the rearranging $4\leftrightarrow 5$, $9\to 7$, $7\to 8$, $8\to 9$. Now we
can apply the triangle-star transformation to the central triangle.}
\label{fig:1.1}
\end{figure}

\begin{figure}
	$$\epsfbox{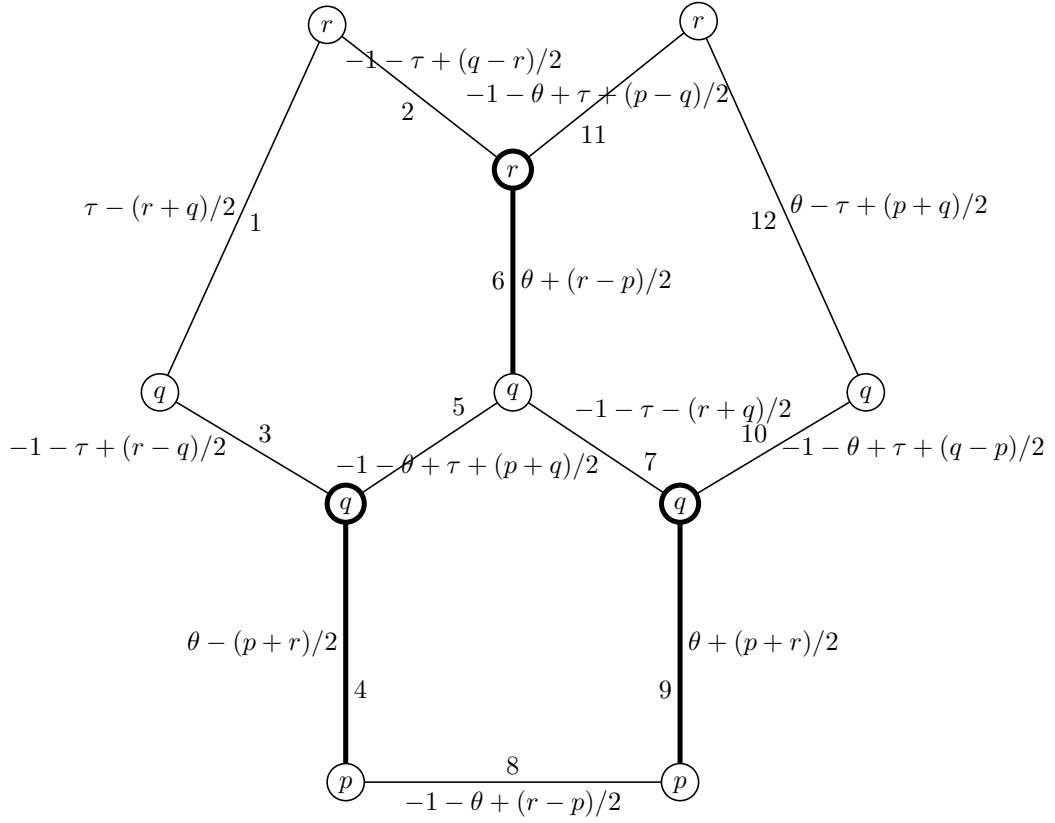}$$
	\caption{After the first star-triangle transformation. We intend to apply star-triangle transformations to emphasized vertices. First, we permute $7\leftrightarrow 8$.}
	\label{fig:2}
\end{figure}

\begin{figure}
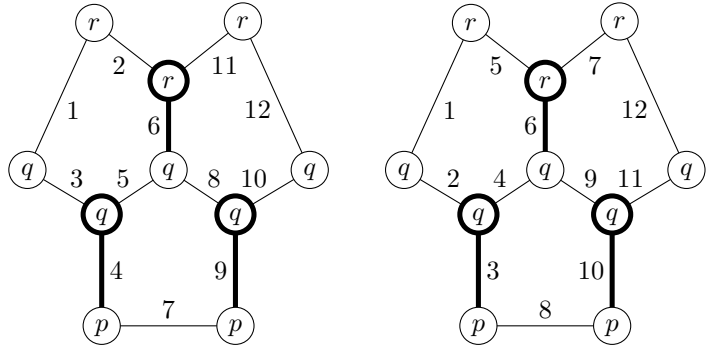

$$\epsfbox{fey7.1}\qquad \epsfbox{fey7.2} $$
\caption{We transpose 2 with product of 3, 4, 5 and 11 with product of 7, 8, 9, 10.
Next, we apply 3 star-triangle transformations simultaneously.}
\end{figure}

\begin{figure}
 $$\epsfbox{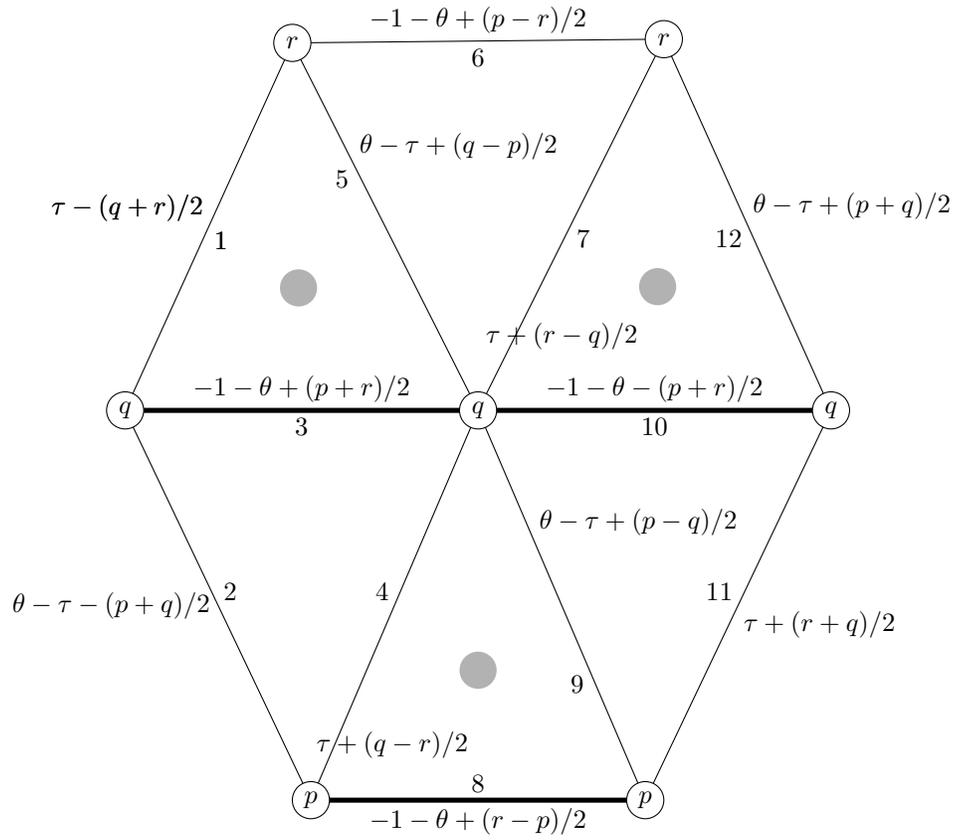}$$
 \caption{After four  transformations. We  transpose 4 with the product of 5, 6, 7.}
 \label{fig:3}
\end{figure}

\begin{figure}
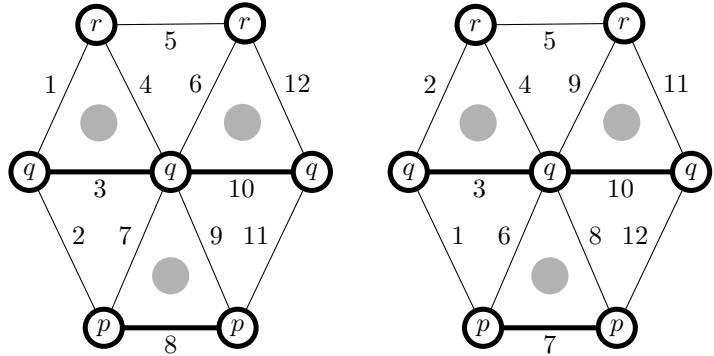

$$\epsfbox{fey8.1}\qquad \epsfbox{fey8.2} $$
\caption{We transpose $1\leftrightarrow 2$, $11\leftrightarrow 12$ and 6 with product of 7, 8, 9.
Now we can do three triangle-star transformations in emphasized triangles.}
\label{fig:predposl}
\end{figure}

\begin{figure}
$$\epsfbox{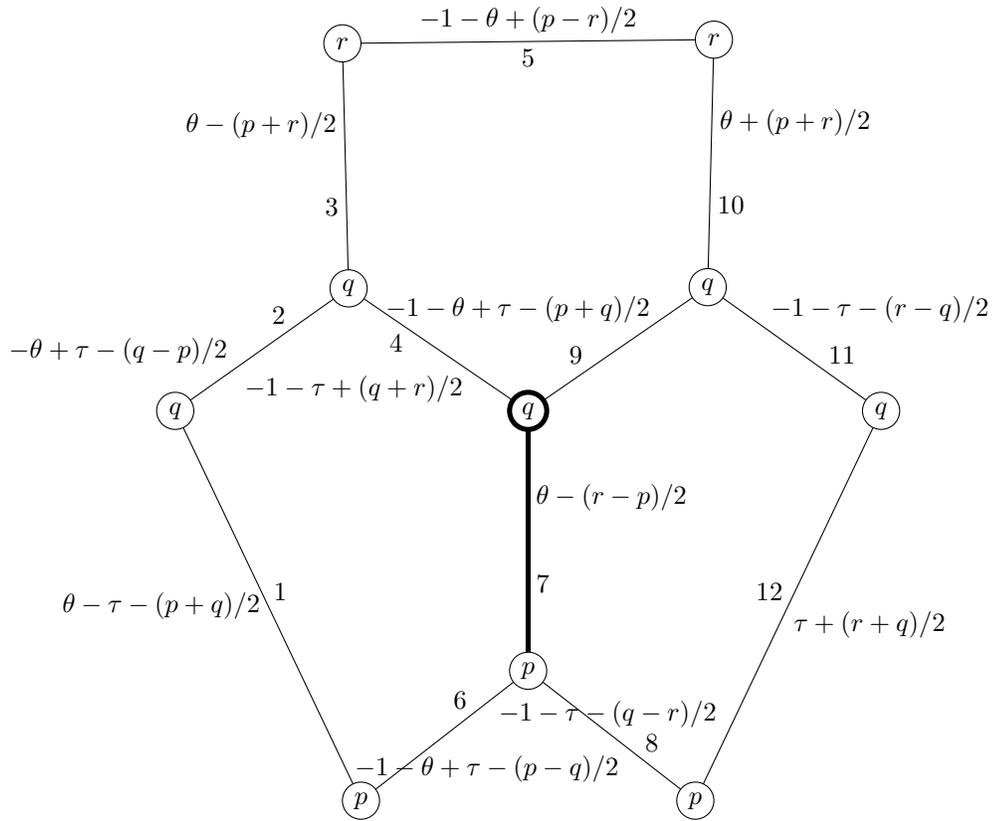}$$
\caption{After 7 transformations. Next, we transpose 4 with product of 5,6 and 
$9\leftrightarrow 8$. After this, we apply the star-triangle transformation and come to the right hand part
of the Yang-Baxter equation \eqref{eq:RRR}.}

\end{figure}

Our expression is a product of 12 operators, formally we can represent \eqref{eq:RRR} as an integral operator
$L^2(\T^3)\to L^2(\T^3)$, whose kernel $K$ is an integral over $\T\times\T\times\T$
and the integrand
is a certain product
of differences of variables in some powers%
\footnote{Recall that the Gauss hypergeometric function $_2F_1$ can be represented as an integral of 
a product $(x-a)^\alpha(x-b)^\beta(x-c)^\gamma$, see, e.g., \cite{HTF}, Sect. 2.1.3, \cite{AAR}, Sect. 2.2. Functions $_pF_{p-1}$ admit integral representations 
as $(p-1)$-dimensional integrals of products of powers of linear forms (we iterate formula (2.2.4) from \cite{AAR}). Some multivariate analogues of hypergeometric functions
(as Appell, Kamp\'et de Feri\'et, Lauricella hypergeometric functions \cite{Ext}, Jack symmetric functions \cite{OO}) 
have integral representations in this spirit.
  So the kernel of the left-hand side of the Yang-Baxter identity can be regarded as a kind of a hypergeometric function of 6 variables,
  and the identity itself is a nontrivial hypergeometric transformation.}.
  We draw a graph, whose vertices correspond to variables (i.e., to coordinates on the initial torus $\T^3$,
 the target torus $\T^3$ and 3 additional copies of $\T$).
Two vertices $v_i$, $v_j$ are connected by an edge iff $K$ contains a factor $|v_i-v_j|^{\mu}$. We write $\mu$ on the corresponding edge, see Fig. \ref{fig:add1}.
 Each edge corresponds to
an operator, products of operators are ordered, and we write on edges numbers of factors (starting from the right). Next,
each variable corresponds to one of the representations $T_p$, $T_q$, $T_r$. So we mark each vertex with a label $p$, $q$, $r$.
In this way, we come to the diagram   Fig. \ref{fig:1}. 

On Fig. \ref{fig:1} we also mark variables corresponding to the initial $L^2(\T^3)$ by entering  arrows,  variables corresponding to the target $L^2(\T^3)$
by outgoing arrows.  Three quadrilaterals correspond to 3 operators $\check R(\cdot)$.

We can transpose numbers $m$ and $m+1$ on edges if the corresponding operators commute. In particular, on Fig. \ref{fig:1} operators  4 and 5 commute, therefore we can transpose
$4\leftrightarrow 5$. The operator  7 commutes with operators 7 and 8,
and  we rearrange  $9\to 7$, $7\to 8$, $8\to 9$, see Fig. \ref{fig:1.1}. Now in the central triangle we get neighboring numeric labels $5$, $6$, $7$, and this allows
us to apply the triangle-star transformation, see Fig.\ref{fig:star-triangle}. We come to the diagram Fig.\ref{fig:2}. Notice that a formal triple integral defining
the kernel was replaced by a quadruple integral.
The whole calculation is drawn on Fig. \ref{fig:1}--\ref{fig:predposl}.


We  can also start with the right hand side of the Yang--Baxter equation. After 4 transformations we come to the hexagon Fig. \ref{fig:3}. But it seems that numeric labels 
require  watching.
\hfill $\square$


{\bf\punct Some variations.} For the group $\SL(2,\C)$ picture is the same. Instead of the integral  \eqref{eq:beta} we must take
the beta-integral \cite{DM}, (2.27). Instead of functions on $\T$ we must consider functions of the two-dimensional sphere (or for non-spherical
representations, 
sections of line bundles on the Riemannian sphere). 
For the spherical principal series of $\SO(1,n)$ we must take the beta-integral 
\cite{CDI}, (5.31).

\section{Problem of definition of operators $\check R(\sigma)$ for complex $\sigma$.}

\COUNTERS

{\bf \punct Analytic continuation.}
Consider a torus $\T^m$. Let $(\alpha_1, \dots, \alpha_n)\in \C^n$. Consider a formal product of operators of the form 
\begin{equation}
V(\alpha_1,\dots, \alpha_n):=B_n \dots B_1,
\label{eq:product}
\end{equation}
where each factor $B_j$ is $A_{pq}(\alpha_j)$ or $J_r(\alpha_j)$. If all $\alpha_j\in\R$, then this operator 
is unitary.
Let $\cL$ be a linear subspace in $\C^n$ determined by a collection of linear equations
with real coefficients.  Let $\cL_{i\R}\subset \cL$ be the subspace consisting of vectors with $\alpha_j\in i\R$.

\begin{theorem} 
\label{th:analytic}
Fix a subspace $\cL$ as above.
 For any $\Phi$, $\Psi\in C^\infty(\T^k)$
 the matrix element
 \begin{equation}
\bigl\la V(\alpha_1,\dots, \alpha_n) \Phi, \Psi\bigr\ra_{L^2(\T^k)}\Bigr|_{\cL_{i\R}}
 \end{equation}
 admits a meromorphic continuation to the subspace $\cL$.
\end{theorem}

{\sc Remark.} It is clear (see the proof below), that $F:=\bigl\la V(\alpha_1,\dots, \alpha_n) \Phi, \Psi\bigr\ra$
is a meromorphic function on $\C^n$. We do not claim that a restriction of $F$
to $\cL$ is well defined. There is a  possibility that $\cL$ is contained in the indeterminacy set of $F$ (i.e., 
an intersection of the hypersurface of poles and the hypersurface of zeros of a meromorphic function $F$), in the proof we 
define a meromorphic function   on such $\cL$ as a certain directional limit.
\hfill $\boxtimes$

\sm

{\sc Example.} The triple product $U_{p,q,r}(\theta,\tau):=\check  R_{23}(\theta-\tau)\check R_{12}(\theta) \check R_{23}(\tau)$ 
in the left hand part of the Yang-Baxter identity \eqref{eq:RRR} is a product
of the form \eqref{eq:product} with 12 factors. So, we have a 5-dimensional subspace $\cL\subset \C^{12}$ parametrized by $p$, $q$, $r$, $\theta$, $\tau$. According the theorem, matrix elements of $U_{p,q,r}(\theta,\tau)$  are well-defined meromorphic functions on
$\cL$. Therefore, we have a meromorphic family of operators $U_{p,q,r}(\theta,\tau)$ acting from $\C^\infty(\T^3)$ to the
space $\cD'(\T^3)$ of all distributions on $\T^3$. For the same reason, we have  three-dimensional meromorphic families 
$\check R_{\dots}^{\dots}(\cdot):\C^\infty(\T^3)\to \cD'(\T^3)$. So all factors $\check R$
and products of these factors
in \eqref{eq:RRR} are well-defined. Now the Yang--Baxter identity as an identity for meromorphic operator-valued functions
seems to be doubtless. 
But {\it I do not understand how to define products of analytic continuations.}
\hfill $\boxtimes$

\sm

{\sc Proof.} 
First, we notice that operator-valued functions $\alpha_j\mapsto A_{pq}(\alpha_j)$ are holomorphic in the domain  $\Re \alpha_j>0$ 
and strongly continuous upto
the boundary $\Re \alpha_j=0$. Operator-valued functions $\alpha_j\mapsto J_{r}(\alpha_j)$ are meromorphic $\Re \alpha_j>0$ and strongly continuous upto
the boundary. Therefore, the function $V(\alpha, \dots,\alpha_n)$ is meromorphic in the domain $\Re \alpha_j>0$ 
and strongly continuous up to the boundary $\Re\alpha_j=0$ (see, e.g, \cite{RS}, Ex. VI.6.c).

Let us write formally the kernel $K(z)$ of the operator $V(\alpha_1,\dots, \alpha_n)$.
It is an integral of a  product of factors $|z_{r_l}-z_{t_l}|^{\mu_l}$, where $\mu_l$ is $\alpha_l$ or $-1+\alpha_l$.
Consider the integral 
\begin{multline}
F(\alpha_1,\dots, \alpha_n) :=\bigl\la V(\alpha_1,\dots, \alpha_n) \Phi,\, \Psi\bigr\ra_{L^2(\T^k)}
=\\=
\int K(z)\, \Phi(z)\, \ov{\Psi(z)}\, \prod \frac{dz_j}{iz_j},
\label{eq:matrix-element}
\end{multline}
where $\Phi(z)$ actually depends on coordinates $z_p$ corresponding to  the first copy of $\T^m$ and $\Psi(z)$
depends on coordinates corresponding to the second copy of $\T^m$. In the domain $\Re \alpha_j>0$  the expression $F(\alpha_1,\dots, \alpha_n)$ is well defined, meromorphic,
and has  the form
\begin{equation}
\gamma(\mu)\,
\int \prod_l |z_{r_l}-z_{t_l}|^{\mu_l}\, \Theta(z)\,\prod \frac{dz_h}{iz_h},
\label{eq:int}
\end{equation}
where $\Theta(z)$ is a smooth function, and $\gamma(\mu)$ is a Gamma-factor.
According  \cite{BG} (see, also, \cite{Sab}, \cite{Gyo})
such expressions have meromorphic
continuations in powers $(\mu_1,\dots, \mu_n)$.
Poles of the integral are located on hyperplanes of
the form 
\begin{equation}
\sum_j b^{(1)}_j \mu_j+c^{(1)}+2s+1=0,\quad\dots,\quad \sum_j b^{(N)}_j \mu_j+c^{(N)}+2s+1=0,
\label{eq:hypeplane}
\end{equation}
where $b^{(l)}_j$, $c^{(l)}$ are fixed non-negative integers, and $s$ ranges in non-negative integers.
Recall that we have $\mu_j=\alpha_j$ or $-1+\alpha_j$, therefore the Gamma-factor (a product of several expressions $C(\alpha_j)^{-1}$, see \eqref{eq:C})
 is holomorphic in a neighborhood of the point $\alpha_1=\dots=\alpha_n=0$.

If $\cL$ is not contained in a hyperplane of type \eqref{eq:hypeplane}, then we simply restrict a meromorphic function to a subspace.
So let  $\cL$ be contained in an intersection of such  hyperplanes. Consider the line $\ell:\alpha_1=\dots=\alpha_n$,
clearly, $\ell\not\subset \cL$. 
Set $\cL^\circ:=\cL\oplus \ell$. We introduce  coordinates $u_1$, $u_2$, \dots on $\cL$ and the coordinate 
$z=\alpha_1$ on $\ell$. 
 The restriction $R(z,u)$ of a function \eqref{eq:matrix-element} to $\cL^\circ$  is well-defined
and meromorphic (since $\cL^\circ$ is not contained in a hyperplane of type \eqref{eq:hypeplane} and 0 is a regular point of the Gamma-factor)
with a possible pole on $\cL$ 
(and possible poles on some other hyperplanes $M_\sigma\subset \cL^\circ$). 
Consider a  point
$\xi_0\in \cL_{i\R}\setminus \cup M_\sigma$. In a neighborhood of $\xi_0$ we have an expansion
$$
R(z,u)=z^{-l} Q_{-l}(u)+ z^{-l+1} Q_{-l+1}(u)+\dots,
$$
where $Q_j(u)$ are holomorphic in the neighborhood  of $\xi_0$. Without loss of generality, we can assume that $l>0$
(otherwise the expression is holomorphic near $\xi_0$). Let $\xi\in \cL_{i\R}$ be in a neighborhood of $\xi_0$.
Then 
$$
\lim_{z\to 0,\, \Re z>0} R(z,\xi)
$$
is finite. Since $R(z,u)$ is meromorphic on the line $\xi+\ell$, its singularity at $z=0$ is removable. So,
$Q_{-l}(\xi)=0$. By the uniqueness theorem, $Q_{-l}(u)=0$ identically.  
Therefore, the pole of $R(z,u)$ on the hyperplane $z=0$ is removable.
\hfill $\square$

\sm

{\bf\punct In which functional spaces can the  operators $\check R(\sigma)$   be defined?}
Theorem \ref{th:analytic} shows that all operators, which are present in formulas \eqref{eq:R}--\eqref{eq:RRR},
are well-defined as operators sending $C^\infty$-functions on a torus to distributions on the same torus.
But this is not sufficient for a definition of products of such operators.

Operators $\check R(\sigma)$ have the form
\begin{equation}
A_{12}(\alpha)\, J_1(\beta)\, J_2(\gamma)\, A_{12}(\delta),\qquad \text{where $\alpha+\beta+\gamma+\delta=0$.}
\label{eq:abcd}
\end{equation}
Operators $J(\beta)$, $J(\gamma)$ are pseudo-differential, such operators have good properties in Sobolev spaces  and other functional spaces
	of Sobolev type, see, e.g., \cite{Tay}, \cite{Shu}. On the other hand, multiplications by non-smooth functions 
	in such spaces is a difficult  topic, which was a subject  of numerous works, see, e.g.,  \cite{MS}, \cite{Mich}, \cite{BS1}. 
	Such operators
	are defined
	in Sobolev spaces only under strong restrictions. As far as I know, existing theorems 
	do not allow us to define products \eqref{eq:abcd}.
	
	However, there is some hope for  the following exotic way  \eqref{eq:abcd}. 
	
	Denote by $|x_1-x_2|^\nu\, C^\infty\,(\T^2)$ the space of functions on torus of the form
	$|x_1-x_2|^\nu f(x_1,x_2)$,  where $\nu\in\C$. We assume that
	the power $\nu$ and all other powers below are $\ne -1$, $-2$, $-3$, \dots.
	
	 We start with a space $|x_1-x_2|^\mu\, C^\infty\,(\T^2)$.
	 
\sm	 
	 
	1) The operator $A_{12}(\delta)$ sends the space $|x_1-x_2|^{\mu}\,C^\infty(\T^2)$
	  to 
	  $$|x_1-x_2|^{\delta+\mu}\, C^\infty\,(\T^2).$$

\sm	 
	 
	 2)	The operator
	$J_2(\gamma)$ sends  the space
$|x_1-x_2|^{\delta+\mu}\, C^\infty\,(\T^2)$ to the space	
$|x_1-x_2|^{\gamma+\delta+\mu}\, C^\infty(\T^2)+C^\infty(\T^2)$,
see, e.g., \cite{Sam}, Sect. 16.


	 3) The operator $J_1(\beta)$  sends the space $|x_1-x_2|^{\gamma+\delta+\mu}\,C^\infty(\T^2)$
	  to the space  $|x_1-x_2|^{\beta+\gamma+\delta+\mu}\, C^\infty(\T^2)+C^\infty(\T^2)$. Also, it sends the space $C^\infty(\T^2)$ to itself.
	 So we get an element of the space
	 	$$|x_1-x_2|^{\beta+\gamma+\delta+\mu}\, C^\infty(\T^2)+C^\infty(\T^2).$$
	
	4) 	Finally, applying $A(\alpha)$ we come to
	 	$$|x_1-x_2|^\mu\, C^\infty(\T^2)+|x_1-x_2|^{\alpha}\, C^\infty(\T^2).$$
	 	 So we can define an operator \eqref{eq:abcd} as an operator 
	 	 $$|x_1-x_2|^\mu\, C^\infty(\T^2) \,
	 	 \longrightarrow \, |x_1-x_2|^\mu\, C^\infty(\T^2)+|x_1-x_2|^{\alpha}\, C^\infty(\T^2).$$
	 	 
	 	 \sm
	 	 
	 	   However, these arguments are not sufficient for a definition of products of  such operators.

\tt

Fakult\"at f\"ur Mathematik, Universit\"at Wien; Institute for Information Transmission Problems;
Faculty of Mechanics and Mathematics, Moscow State University.

yurii.neretin(frog)univie.ac

URL: https://mat.univie.ac.at/$\sim$neretin/

\end{document}